\documentclass[12pt]{article}
\usepackage{makeidx}
\usepackage{amssymb}
\usepackage{amsfonts}
\usepackage{amsmath}
\usepackage{euscript}
\usepackage{theorem}
\usepackage{enumerate}

\begingroup

\newtheorem{thm}{Theorem\hskip 5mm}[section]

\newtheorem{prop}[thm]{Proposition\hskip 5mm}

\newtheorem{cor}[thm]{Corollary\hskip 5mm}

\newtheorem{lem}[thm]{Lemma\hskip 5mm}

\newtheorem{exa}[thm]{Example\hskip 5mm}

\newtheorem{claim}[thm]{Claim\hskip 5mm}

\endgroup

\begingroup

\endgroup

\begingroup

\endgroup

\def\go{{\sqrt{5}}}

\def\GL{{\mathrm {GL}}}

\def\Q{{\bf Q}}
\def\C{{\bf C}}
\def\Z{{\bf Z}}

\def\a{{\alpha}}
\def\b{{\beta}}

\begin{document}

\begin{center}
{\bf IS EVERY MATRIX SIMILAR TO A POLYNOMIAL\\ IN A
COMPANION MATRIX?}
\end{center}

\begin{center}{{\large {\sc N. H. Guersenzvaig}}\\
{\small A\!v. Corrientes 3985 6A, (1194) Buenos Aires, Argentina\\ email: nguersenz@fibertel.com.ar}}
\end{center}

\begin{center}
and
\end{center}

\begin{center}{{\large {\sc Fernando Szechtman}}\\
{\small Department of Mathematics and Statistics, University of Regina, Saskatchewan, Canada\\
email: fernando.szechtman@gmail.com}}
\end{center}

\begin{abstract}

Given a field $F$, an integer $n\geq 1$, and a matrix
$A\in M_n(F)$, are there polynomials $f,g\in F[X]$, with $f$ monic
of degree $n$, such that $A$ is similar to $g(C_f)$, where $C_f$
is the companion matrix of $f$? For infinite fields the answer is easily seen to positive,
so we concentrate on finite fields. In this case we give an affirmative answer, provided
$|F|\geq n-2$. Moreover, for any finite field $F$, with $|F|=m$,
we construct a matrix $A\in M_{m+3}(F)$ that is not similar to any
matrix of the form $g(C_f)$.

Of use above, but also of independent interest, is a constructive procedure
to determine the similarity type of any given matrix
$g(C_f)$ purely in terms of $f$ and $g$, without resorting to
polynomial roots in $F$ or in any extension thereof. This, in turn, yields
an algorithm that, given $g$ and the invariant factors of any $A$,
returns the elementary divisors of~$g(A)$. It is
a rational procedure, as opposed to the classical method that uses the Jordan
decomposition of $A$ to find that of~$g(A)$.


Finally, extending prior results by the authors, we show that for an
integrally closed ring $R$ with field of fractions $F$ and
companion matrices $C,D$ the subalgebra $R\langle C,D\rangle$ of
$M_n(R)$ is a free $R$-module of rank $n+(n-m)(n-1)$, where $m$ is
the degree of $\gcd (f,g)\in F[X]$, and a presentation for
$R\langle C,D\rangle$ is given in terms of $C$ and~$D$. A
counterexample is furnished to show that $R\langle C,D\rangle$
need not be a free $R$-module if $R$ is not integrally closed.
The preceding information is used to study $M_n(R)$, and others,
as $R[X]$-modules.
\end{abstract}

\medskip


\smallskip


\smallskip

\small{The second author was supported in part by an NSERC
discovery grant}

\section{Introduction}

We fix throughout the paper a field $F$ and an integer $n\geq 1$.
Given a monic polynomial $f=f_0+f_1X+\cdots+f_{n-1}X^{n-1}+X^n\in
F[X]$, its companion matrix $C_f\in M_n(F)$ is defined:
$$
C_f=\left(%
\begin{array}{ccccc}
  0 & 0 & \cdots & 0 & -f_0 \\
  1 & 0 & \cdots & 0 & -f_1 \\
  0 & 1 & \cdots & 0 & -f_2 \\
  \vdots & \vdots & \cdots & \vdots & \vdots \\
  0 & 0 & \cdots & 1 & -f_{n-1} \\
\end{array}%
\right).
$$

Companion matrices play a prominent role in linear algebra.
Indeed, any matrix $A\in M_n(F)$ is similar to a unique direct sum
$C_{q_1}\oplus\cdots\oplus C_{q_r}$, where $1\neq q_1|\cdots|q_r$
are the invariant factors of~$A$. Refining this decomposition
results in $A$ being similar to $A_{p_1}\oplus\cdots\oplus
A_{p_s}$, where $P=\{p_1,\dots,p_s\}$ is the set of monic
irreducible factors of the minimal polynomial of $A$ and, for
$p\in P$, the component $A_p$ is similar to
$C_{p^{a_1}}\oplus\cdots\oplus C_{p^{a_t}}$ for unique $1\leq
a_1\dots\leq a_t$. The powers $p^{a_1},\dots,p^{a_t}$ are known as
the $p$-elementary divisors of $A$, and as $p$ runs through $P$
one obtains all elementary divisors of $A$. This material is
classical and can be found, for instance, in \cite{J}.

Companion matrices and their applications are featured extensively
in the mathematical literature, not only in linear algebra, but
also elsewhere. The list is too long to describe here, so we
restrict ourselves to a sample of cases.

In linear algebra itself, \cite{DW} gives a polar decomposition
for $C_f$, \cite{KS} provides a singular value decomposition as
well as the $QR$-factorization for $C_f$, while the the numerical
range of $C_f$ is studied in \cite{GW}. A new list of
representatives under matrix similarity was recently given in
\cite{D} using companion matrices.

Since $C_f$ is the ``universal root" of $f$, companion matrices
are naturally associated to simple algebraic field extensions.
Not long ago, \cite{GS} gave an application in this
direction. Let $\a\in\C$ be an algebraic number whose minimal
polynomial over $\Q$ is assumed to be $f$. Then
$1,\a,\dots,\a^{n-1}$ is a $\Q$-basis of the field $\Q[\a]$. If
$\b,\gamma\in\Q[\a]$ then \cite{GS} gives a closed formula for the
coordinates of $\b\gamma$ in terms of the individual coordinates
of $\b$ and $\gamma$ and the companion matrix $C_f$. The formula
extends beyond the field context to arbitrary commutative rings
with identity.

For applications of companion matrices to the solution of
differential equations see \cite{B2}, to dynamical systems see
\cite{LS}, to quasicrystals see \cite{GV}, to modular
representation theory see \cite{N}, to the irreducibility of
polynomials see \cite{Gu}, to the study of zeroes of polynomials
see \cite{CN}.

Given monic polynomials $f,g$ of degree $n$ over a commutative
ring with identity $R$, let $C$ and $D$ stand for their
respective companion matrices. Recently, \cite{GS2} studied the
subalgebra $R\langle C,D\rangle$ of $M_n(R)$ generated by $C$ and
$D$. It is shown in \cite{GS2} that $R\langle C,D\rangle$ always
coincides with $R[C,D]$, the $R$-span of all matrices $C^iD^j$,
where $0\leq i,j<n$. Moreover, necessary and sufficient conditions
for $C,D$ to generate $M_n(R)$ were given, and a presentation for
$M_n(R)$ in terms of $C$ and $D$ was furnished in that case.
Furthermore, if $R=\Z$ and $f,g$ are relatively prime then the
sublattice $\Z\langle C,D\rangle$ of $M_n(\Z)$ was shown to have
full rank $n^2$ and the exact value of the finite index
$[M_n(\Z):\Z\langle C,D\rangle]$ was found to be $|R(f,g)|^{n-1}$,
where $R(f,g)$ is the resultant to $f$ and $g$.

A matrix $A\in M_n(F)$ is said to be cyclic if it is similar to a
companion matrix, i.e., if the minimal and characteristic
polynomials of $A$ coincide. The terminology arises from the fact that
when one views the column space $V=F^n$ as an $F[X]$-module via
$A$, this module has a cyclic generator if and only if the matrix
$A$ is cyclic. We will say that $A\in M_n(F)$ is of {\em polynomial type} if
$A$ is equal to a polynomial in a
cyclic matrix. Since the only matrices commuting with a cyclic matrix are
are its own polynomials, being of polynomial type means commuting with a cyclic
matrix. Alternatively, $A\in M_n(F)$ is of polynomial type if and only if $A$ is similar to a matrix $g(C_f)$ for suitable
polynomials $f,g\in F[X]$, with $f$ monic of degree $n$.

In 1995 Neumann and Praeger \cite{NP} studied the
proportion of cyclic matrices in $M_n(F)$ when $F$ is a finite
field. This initiated a long series of papers (see \cite{W},
\cite{NP2}, \cite{Br}, \cite{FNP}, for instance), that study the
probability of a matrix $A$ to be cyclic when it belongs to a
classical linear group defined over a finite field. Now, it is shown in \cite{NP} that the probability that a matrix
be cyclic is high. In this paper we consider the following
question: is every matrix $A\in M_n(F)$ of polynomial type? The answer is easily seen to be positive
for infinite fields, so we concentrate on the case when $F$ is finite. In this case Theorem
\ref{rojo} gives an affirmative answer when $|F|\geq n-2$. This is
actually a consequence of the more general Theorem \ref{main},
which gives sufficient conditions for a matrix $A$ of arbitrary
size $n$ over $F$ to be of polynomial type, depending on the
nature of the elementary divisors of $A$ and the size of $F$. For
a while it was unclear whether the bound $|F|\geq n-2$ was the
result of the inadequacy of our methods or a counterexample really
existed. Eventually the mystery was dispelled: this bound is as
sharp as possible. Indeed, for every finite field~$F$, say of size
$|F|=m$, we construct a matrix $A\in M_{m+3}(F)$ that is not of
polynomial type. Many more examples, of both positive
and negative nature, are given throughout the paper.

A basic tool used to deal with matrices of polynomial type is
given in section \ref{jam}, which is of independent interest. We
furnish a constructive procedure to determine the similarity type
of $g(C_f)$ purely in terms of $f$ and $g$. No appeal is made to
roots of $f$ in $F$ or in any extension
thereof. This naturally extends to an algorithm that yields the elementary divisors of any $g(A)$
from $g$ and the invariant factors of $A$, without explicitly computing $g(A)$. It is
a rational procedure, as opposed to the classical method that uses the Jordan
decomposition of $A$ to find that of $g(A)$, as found in \cite{M}. As a corollary we obtain constructive
criteria for $g(A)$ to be semisimple, cyclic or diagonalizable.

At the end of the paper we direct attention to other results from
\cite{GS2}. If $R$ is a unique factorization domain then
\cite{GS2} shows that $R\langle C,D\rangle$ is necessarily  a free
$R$-module of rank $n+(n-m)(n-1)$, and gives a presentation for
$R\langle C,D\rangle$ in terms of $C,D$. Here we push these
results under the much weaker hypothesis that $R$ be integrally
closed. We also show that if $R=\Z[\sqrt{5}]$ then $R\langle C,D\rangle$ is not a free
$R$-module for suitable choices of $f$ and $g$.

\section{Similarity type of $g(A)$}
\label{jam}

Given $g\in F[X]$ and $A\in M_n(F)$ we wish to determine the elementary divisors of $g(A)$
without computing $g(A)$, finding roots or leaving the field $F$ in any way.

We start by making some observations about the determination of the similarity class of
an arbitrary matrix $A\in M_n(F)$. The interested reader may consult [J] for further information.

We can make the column space $V=F^n$ into an $F[X]$-module via $A$ as follows:
$$
y\cdot v=y(A)v,\quad y\in F[X],v\in V.
$$
For $y\in F[X]$ we define
$$
V(y,A)=\{v\in V\,|\, y\cdot v=0\}=\{v\in V\,|\, y(A)v=0\}.
$$
Suppose $p\in F[X]$ is irreducible of degree $k$. Then $K=F[X]/(p)$ is a field and an $F$-vector space of dimension $k$.
Moreover, $V(p,A)$ is a $K$-vector space and so is $V(p^{i+1},A)/V(p^{i},A)$ for every $i\geq 1$. It follows that
the dimension of $V(p^i,A)$ over $F$ is a multiple of $k$ for every $i\geq 1$. Here $V(p^0,A)=\{0\}$ and
$V(p,A)\neq \{0\}$ if and only if $p$ is a factor of the minimal polynomial $q$ of $A$. Let $q$ have prime factorization $q=p_1^{e_1}\cdots p_m^{e_m}$,
where $P=\{p_1,\dots,p_m\}$ is the set of of monic irreducible factors of $q$ and $e_1,\dots,e_m\geq 1$.
Let $p\in P$ have multiplicity $e$ in $q$, and define the sequence $d_0,d_1,d_2,\dots$ by
$$
d_i=\mathrm{dim}_F(V(p^i,A))/\mathrm{deg}(p),\quad i\geq 0.
$$
Here $d_0=0$ and each $d_i$, for $i\geq 1$, is a positive integer, as explained above. Note also that $V(p^e)=V(p^{e+1})=\cdots$, so the sequence $d_0,d_1,d_2,\dots$ stabilizes at $d_e$.

\begin{claim}
\label{calmp}
The similarity type of $A$ is completely determined by the $m$ sequences $d_1,d_2,\dots$ corresponding to each of the polynomials $p_1,\dots,p_m$.
\end{claim}

\noindent{\it Proof.} Let $p\in P$ have multiplicity $e$ in $q$. For each $i\geq 1$ let $a_i$ stand for the number of elementary divisors of $A$ of the form $p^k$ with $k\geq i$ and let $b_i$ stand for the number of the form~$p^i$. Clearly $b_1=a_1-a_2, b_2=a_2-a_3, b_3=a_3-a_4\dots$. On the other hand, looking at the decomposition
of $A$ into elementary divisors, we see that  $$d_1=a_1, d_2-d_1=a_2, d_3-d_2=a_3, d_4-d_3=a_4\dots,$$ so after all
$$
b_1=2d_1-d_2-d_0, b_2=2d_2-d_3-d_1, b_3=2d_3-d_4-d_2,\dots
$$
i.e.,
\begin{equation}
\label{dbeq}
b_i=2d_i-d_{i+1}-d_{i-1},\quad i\geq 1.
\end{equation}
Since the multiplicities of the elementary divisors of $A$ completely describe $A$ up to similarity, our claim
is established.$\quad
\blacksquare$

Note that since all terms from $d_e$ on are equal, it follows that $b_i=0$ for $i\geq e+1$, as expected.

Next let $f\in F[X]$, where $f$ is monic of degree $n$. We wish to apply the above ideas to determine the
elementary divisors of $g(C_f)$ solely in terms of $f$ and $g$. We assume that we have somehow
found the set $R=\{r_1,\dots,r_k\}$ of monic irreducible factors of $f$.
A polynomial time algorithm for factoring polynomials over the rational numbers
can be found in \cite{LLL}, and over an algebraic number field in \cite{L}. Factorization methods
for polynomials over a finite field can be found in \cite{Be} and \cite{CZ}. Let $r\in R$ have degree $s$. By the minimal polynomial of $g$ modulo~$r$ we  mean the monic generator $p$
of the ideal of all $y\in F[X]$ such that $y(g(X))$ is a multiple of $r$. To calculate $p$
simply reduce the powers $1,g(X),g(X)^2,\dots$ modulo~$r$. Then
the first power $g(X)^t$, with $t\leq s$, that is a linear combination, modulo $r$, of those preceding it will
yield the actual coefficients to construct $p$. Notice that $\mathrm{deg}(p)\leq \mathrm{deg}(r)$
and that $p$ is irreducible since so is $r$. It is perfectly possible for $g(X)$ to have
the same minimal polynomial modulo distinct factors $r_i$ and $r_j$ of $f$.

We claim that the monic irreducible factors of the minimal polynomial $q$
of $g(C_f)$ are precisely the minimal polynomials of $g$ modulo
the irreducible factors of $f$. Indeed, for $y\in F[X]$ we have
$y(g(C_f))=0$ if and only if $y(g(X))\equiv 0\mod f$, which is
equivalent to $y(g(X))\equiv 0\mod r_i^{u_i}$, $1\leq i\leq k$, where
$u_i$ is the multiplicity of $r_i$ in $f$. This readily yields our
claim.

Having determined the monic irreducible factors $p$ of $q$
what we need now is a way to compute the dimension of
$V(p^i,g(C_f))$ for $i\geq 1$. We will actually compute $\mathrm{dim}\; V(y,g(C_f))$ for every $y\in F[X]$. This requires some notation.

Let $e_1, \dots , e_n$ stand for the canonical basis of $V$. We write $F_n[X]$ for the subspace of $F[X]$
with basis $1,X,\dots,X^{n-1}$. If $p\in F_n[X]$ then $[p]$ stands
for the coordinates of $p$ relative to this basis. Note that every vector $v\in V$ can be
written in the form $v=[p]$ for a unique polynomial $p\in F_n[X]$.

\begin{lem}
\label{uno} Given any $y\in F[X]$, we set $z=\gcd(y(g(X)),f)$, $d=\mathrm{deg}(z)$, $h=f/z$, $V(y)=V(y,g(C_f))$
and $C=C_f$. Then

(a) If $p\in F_n[X]$ then $[p]$ is in $V(y)$ if and only if $h|p$.

(b) The dimension of $V(y)$ is equal to $d$, i.e.,
$$
\mathrm{deg}\; gcd(y(g(X)),f)=\mathrm{dim}\; V(y,g(C_f)).
$$
\indent (c)  If $d>0$ then a basis for $V(y)$ is given by $[h],[X h],\dots, [X^{d-1} h]$.
\end{lem}

\noindent{\it Proof.} By definition $[p]\in V(y)$ if and only if
$$
y(g(C))[p]=0.
$$
We easily verify that $[p]=p(C)e_1$, so the above translates into
\begin{equation}
\label{toto}
y(g(C))p(C)e_1=0.
\end{equation}
The polynomials that evaluated in $C$ annihilate $e_1$ are those that annihilate $C$, namely the multiples of $f$.
Thus (\ref{toto}) is equivalent to $f|y(g(X))p(X)$, i.e. $h|p$. This proves (a).
It is clear then that $[h],[X h],\dots, [X^{d-1} h]$ are all in $V(y)$. Moreover, looking at the last 1 in each
of these column vectors, we see that they are linearly independent. The conditions $h|p$ and $p\in F_n[X]$ show that
$[h],[X h],\dots, [X^{d-1} h]$ span $V(y)$, as required.$\quad
\blacksquare$

By Lemma \ref{uno}, if $p$ is one of the irreducible factors of the minimal polynomial of $g(C_f)$ then
$$
\mathrm{deg}\; \gcd(p^i(g(X)),f)=\mathrm{dim}\; V(p^i,g(C_f)),\quad i\geq 1.
$$

By our preceding discussion the values on the right hand side completely determine the $p$-elementary divisors of $g(C_f)$ and hence so do those on the left. From a computational point of view, note that these degrees will strictly increase until a first unavoidable repetition; all values onward will coincide, with no need to compute them.
We have thus proved the following.

\begin{thm}
\label{k}
The similarity type of $g(C_f)$ is completely determined by the set $P$ of minimal polynomials $p$ of $g$ modulo
the monic irreducible factors $r$ of $f$ and the degrees of
$$\gcd(p(g(X)),f),\;\gcd(p^2(g(X)),f),\;\gcd(p^3(g(X)),f),\dots
$$
Explicitly, the set of monic irreducible factors of the minimal polynomial of $g(C_f)$ is $P$ and for $p\in P$, if
$d_i=\mathrm{deg}\;\gcd(p^i(g(X)),f)/\mathrm{deg}\; p$, $i\geq 1$, then the multiplicity $b_i$ of $p^i$ as an elementary divisor of $g(C_f)$
is given by (\ref{dbeq}).
\end{thm}
Observe next that if $y\in F[X]$ is an arbitrary polynomial and
$h=\gcd(y,f)$ then
\begin{equation}
\label{pat} \gcd (y^i,f)=\gcd (h^i,f),\quad i\geq 1.
\end{equation}
This can be used to our advantage to ease the computation of $\gcd(p^i(g(X)),f)$ for $i>1$.

It is now easy to complete the determination of the similarity type of $g(A)$ from that of $A$.
We have $A\sim C_{f_1}\oplus\cdots\oplus C_{f_t}$, where $1\neq f_1|\cdots|f_t$ are
the invariant factors of $A$ and $f_t$ is the minimal polynomial of $A$. Clearly
the minimal polynomial of $g(A)$ is the minimal polynomial $q_t$ of $g(C_{f_t})$.
We compute all monic irreducible factors $p$ of $q_t$ and, for each of these,
the corresponding $p$-elementary divisors of $g(C_{f_1}),\dots,g(C_{f_t})$,
as indicated above. This produces the $p$-elementary divisors of $g(A)$, as required.

\begin{prop}
\label{z2} Suppose $g(X)$ has degree $d\geq 1$
and leading coefficient $a$. Then
$$
g(C_{a^{-n} f(g(X))})\sim C_f\oplus\cdots\oplus C_f,\quad d\text{ times }
$$
\end{prop}

\noindent{\it Proof.} Let $r$ be a monic irreducible factor of $f$
with multiplicity $e$, and let $z$ be a monic irreducible factor
of $r(g(X))$. The minimal polynomial of $g(X)$ modulo $z$ is
clearly just $r$. As $r$ varies in $f$ and $z$ varies in
$r(g(X))$, the resulting $z$ exhaust all monic irreducible factors
of $a^{-n} f(g(X))$. Thus, the monic irreducible factors of the
minimal polynomial of $g(C_{a^{-n}f(g(X))})$ are precisely those
of $f$. Moreover, if $1\leq i\leq e$ we have
$$
\mathrm{deg}\; \gcd(r^i(g(X)),f(g(X)))=\mathrm{deg}\; r^i(g(X))=d\times \mathrm{deg}(r^i)=d\times \mathrm{deg}(\gcd(r^i,f)),
$$
while if $i>e$ we have
$$
\mathrm{deg}\;\gcd(r^i(g(X)),f(g(X)))=d\times \mathrm{deg}(r^e)=d\times \mathrm{deg}(\gcd(r^i,f)).
$$
The result now follows from Theorem \ref{k}.$\quad
\blacksquare$

\section{Matrices of polynomial type}
\label{det}

Recall the following well-known result.

\begin{lem}
\label{comiu}
 Suppose $A\in M_n(F)$ is cyclic and $B$ commutes with $A$. Then $B\in F[A]$.
\end{lem}

If $A,B\in M_n(F)$ we write $A\sim B$ to denote that $A,B$ are similar. If $F$ is an algebraically closed field, the fact
that every $A\in M_n(F)$ is of polynomial type has been known for a long time. See, for instance, the argument given in \cite{MT}, page 111.
The proof below, due to R. Guralnick, appears in \cite{GR}.

\begin{thm}
\label{thor}
Let $F$ be field. Suppose all eigenvalues of $A\in M_n(F)$ are in $F$,
and $F$ has at least as many elements as the number of Jordan blocks of $A$. Then
$A$ is of polynomial type.
\end{thm}

\noindent{\it Proof.} By hypothesis $A\sim B=J_{m_1}(a_1)\oplus\cdots\oplus J_{m_s}(a_s)$, where
$a_1,\dots,a_s\in F$. By assumption
$F$ has at least $s$ distinct elements $b_1,...,b_s$.
Let $D=J_{m_1}(b_1)\oplus\cdots\oplus J_{m_s}(b_s)$.
Then $B$ commutes with~$D$, which is cyclic, so $B\in F[D]$ by Lemma \ref{comiu}, whence $A$ is of polynomial type.
$\quad\blacksquare$

\begin{cor}
\label{thor2}
If $F$ is algebraically closed then every $A\in M_n(F)$ is of polynomial type.
\end{cor}

The general case is not as easy. We proceed to tackle this problem.

\begin{lem}
\label{chin}
 Suppose $A_1,\dots,A_m\in M_n(F)$ satisfy $A_i\sim g_i(C_{f_i})$, where $g_1,\dots,g_m\in F[X]$, and $f_1,\dots,f_m\in F[X]$
are pairwise relatively prime. Then $A=A_1\oplus\cdots\oplus A_m$ is similar to $g(C_f)$ for some $g\in F[X]$ and
$f=f_1\cdots f_n$.
\end{lem}

\noindent{\it Proof.} By the Chinese Remainder Theorem there exists $g\in F[X]$ such that
$g\equiv g_i\mod f_i$ for all $1\leq i\leq m$. It follows that
$$
g(C_{f_1}\oplus\cdots\oplus C_{f_m})= g_1(C_{f_1})\oplus\cdots\oplus g_m(C_{f_m})\sim A_1\oplus\cdots\oplus A_m= A.
$$
Since $f_1,\dots,f_m$ are pairwise relatively prime, their product $f$ satisfies
$$
C_{f_1}\oplus\cdots\oplus C_{f_m}\sim C_f.
$$
Therefore
$$
g(C_f)\sim g(C_{f_1}\oplus\cdots\oplus C_{f_m})\sim A.\quad
\blacksquare$$

\medskip

We consider now the subgroup $\Omega$ of $\GL_2(F)$ of all matrices of the form $M=\left(\begin{array}{cc} a & b\\ 0 & 1\end{array}\right)$,
where $a,b\in F$ and $a\neq 0$. Let $F^+$ and $F^*$ stand for the underlying additive and multiplicative groups of $F$, respectively.
We can view $F^+$ and $F^*$ as subgroups of $\Omega$, and in fact $\Omega\cong F^+\rtimes F^*$.

The group $\Omega$ acts on the polynomial algebra $F[X]$ by means of automorphisms as follows. If $f\neq 0$ then
$$
(f^M)(X)=a^{-m}f(aX+b),\quad m=\mathrm{deg}(f),
$$
while $0^M=0$. Note that $\Omega$ preserves the class of monic
irreducible polynomials. We will refer to $f^M  $ as an
$\Omega$-conjugate of $f$. We will also speak of the $F^*$ and
$F^+$-conjugates of $f$. We write $S_f$ and $T_f$ for the
stabilizers of $f$ in $F^+$ and $F^*$, respectively.

Suppose $f$ has degree $m\geq 1$ and let $\a$ be a root of $f$ in
some extension of $F$. If $b\in F^+$ then $\a-b$ is a root of
$f^b$. Thus if $b\in S_f$ it follows that $\a-b$ is also a root
of $f$. Clearly this can happen for at most $m$ such $b$, i.e.
$|S_f|\leq \mathrm{deg}(f)$. In particular, $S_f$ is finite.
Assume next that $F$ is a finite field with $q=p^d$ elements, with $p$ a prime, and
that $f$ is irreducible. The roots of $f$ are distinct from each
other and are permuted by $S_f$, although not necessarily in a
transitive way. We claim that at least one of the stabilizers
$S_f$, $T_f$ is trivial. Indeed, let $K=F[\a]$, a field extension
of $F$ of degree $m$. The Galois group $G$ of $K/F$ is cyclic of
order $m$ generated by the Frobenius automorphism $g$ given by
$x\mapsto x^q$. The group $G$ acts transitively on the roots of
$f$. Suppose $S_f$ is not trivial and let $0\neq b\in S_f$. Then
$m\geq 2$ and some power of $g$ must send $\a$ to $\a-b$, i.e.
$\a^{q^e}=\a-b$ for some $1\leq e\leq m-1$. Since $X^{q^e}-X+b\in
F[X]$ annihilates $\a$, it follows that $f|X^{q^e}-X+b$. Let $a\in
T_f$. We wish to show that $a=1$. Since $a^{-1}\a$ is a root of
$f^a$, it is also a root of $f$, and hence of $X^{q^e}-X+b$.
Thus
$$0=(a^{-1}\a)^{q^e}-a^{-1}\a+b=a^{-1}(\a^{q^e}-\a+ab)=a^{-1}(-b+ab).
$$
From $b=ab$ and $b\neq 0$ we infer $a=1$, as claimed. It follows
that the number of $\Omega$-conjugates of $f$ is at least $|F^*|$.

Note that while $X^{q^e}-X+b$ is $F^+$-invariant, its
irreducible factors need not be: they are permuted by the elements of $F^+$ (e.g. $X^2+X+\a\in F_4[X]$, where $\a^2=\a+1$, is irreducible and $F_2$-stable but not $F_4$-stable, as is the product $(X^2+X+\a)(X^2+X+\a+1)=X^4+X+1$).

For instance, $X^p-X+b\in F_p[X]$, $b\neq 0$, is irreducible,
$F_p^+$-invariant and has itself the required form. These are the
only irreducible $F_p^+$-invariant polynomials of prime degree $p$. Next let $F=F_2$. Clearly $X^8+X+1$ is
$F_2^+$-invariant with prime decomposition
$X^8+X+1=(X^2+X+1)(X^6+X^5+X^3+X^2+1)$. Thus the irreducible
polynomial $X^6+X^5+X^3+X^2+1$ is $F_2^+$-invariant and divides
the polynomial $X^8-X+1$, which has the anticipated form.

\begin{thm}
\label{main}
Let $F$ be an arbitrary field and let $A\in M_n(F)$. Given a monic irreducible factor $p$ of the
minimal polynomial of $A$, let $\ell(p)$ be the sum of
the multiplicities of all elementary divisors of $A$ which are a
power of an $\Omega$-conjugate of $p$.

Suppose that $|F^+|\geq \ell(p)$ if $S_p=\{0\}$, and $|F^*|\geq \ell(p)$ if $S_p\neq \{0\}$ but $T_p=\{1\}$, for every monic irreducible factor $p$ of the minimal polynomial of $A$. Then $A$ is of polynomial type.
More precisely, let $q_1,\dots,q_m$ be the elementary divisors of $A$, where
each $q_i=p_i^{e_i}$ is a positive power of a monic irreducible
factor $p_i$ of the minimal polynomial of $A$ (we do not assume at
all here that $p_1,\dots,p_m$ are distinct). Then $A\sim g(C_f)$, for some $g\in F[X]$ and $f=f_1\cdots f_m$,
where each $f_i=r_i^{e_i}$ for a suitable $\Omega$-conjugate $r_i$ of $p_i$, and $r_1,\dots,r_m$ are distinct.
\end{thm}

\noindent{\it Proof.} The elementary divisors of $A$ give a
decomposition $A\sim C_{q_1}\oplus\cdots\oplus C_{q_m}$, where
each $q_i=p_i^{e_i}$, as indicated above. The result is obvious
if $m=1$ so we suppose $m\geq 2$.

Set $r_1=p_1$. Suppose $1\leq s<m$ and we have found
{\em distinct} polynomials $r_1,\dots,r_s$ such that each $r_i$ is $\Omega$-conjugate of $p_i$. We wish to find $r_{s+1}$,
distinct from all $r_1,\dots,r_s$ and an $\Omega$-conjugate to $p_{s+1}$.

Suppose first $F$ is infinite. As explained above, the $F^+$-stabilizer of $p_{s+1}$ is finite, so
$p_{s+1}$ has infinitely many distinct $F^+$-conjugates. Hence there exist infinitely many $\Omega$-conjugates of $p_{s+1}$
which are different from the given polynomials $r_1,\dots,r_s$ and we just pick one of them to be $r_{s+1}$.

Suppose next $F$ is finite. As shown above, one of the stabilizers $S_{p_{s+1}}$, $T_{p_{s+1}}$ must be trivial.
If $S_{p_{s+1}}=\{0\}$ then the number of $F^+$-conjugates of $p_{s+1}$ is $|F^+|$. In the list $1,\dots,s$
there are less than $\ell(p_{s+1})$ indices $i$ such that $r_i$ is $\Omega$-conjugate to $p_{s+1}$ (since $q_{s+1}$,
which contributes to $\ell(p_{s+1})$, has not been counted yet). By assumption $|F^+|\geq \ell(p_{s+1})$,
so we can find an $F^+$-conjugate $r_{s+1}$ to $p_{s+1}$ different from all $r_1,\dots,r_s$. If $S_{p_{s+1}}\neq \{0\}$
then $T_{p_{s+1}}=\{1\}$ and we may argue exactly as above, using $F^*$ instead of $F^+$.

We thus produce  distinct polynomials $r_1,\dots,r_m$ such that
each $r_i$ is $\Omega$-conjugate to $p_i$. Since $\Omega$
preserves the class of monic irreducible polynomials,
$r_1,\dots,r_m$ are pairwise relatively prime, and hence so are
their powers $f_i=r_i^{e_{i}}$, $1\leq i\leq m$. Clearly $f_i$ is
$\Omega$-conjugate to $q_i$ for every $1\leq i\leq m$ (by an
element that is either in $F^+$ or in $F^*$). It follows from
Proposition \ref{z2} that $C_{q_i}\sim g_i(C_{f_i})$ for a suitable
linear polynomial $g_i$. Now apply Lemma \ref{chin}.
$\quad\blacksquare$

\begin{lem}
\label{xr} If $A$ is diagonalizable then $A\sim g(C_f)$ for
suitable $f,g\in F[X]$.
\end{lem}

\noindent{\it Proof.} We have $A\sim a_1I_{m_1}\oplus\cdots\oplus a_k
I_{m_k}$, where $a_1,\dots,a_k$ are distinct elements of $F$. Then
$a_i I_{m_i}=g_i(C_{f_i})$ for $f_i=(X-a_i)^{m_i}$ and $g_i=a_i$,
so Lemma \ref{chin} applies.  $\quad\blacksquare$

\begin{lem}
\label{xr2} If the elementary divisors of $A\in M_n(F)$ are $(X-a_1)^2,(X-a_2),\dots,(X-a_{n-1})$
for some elements $a_1,\dots,a_{n-1}\in F$, not necessarily distinct, then $A$ is of polynomial type.
\end{lem}

\noindent{\it Proof.} We have $A\sim C_{(X-a_1)^2}\oplus a_1 I_{m}\oplus b_1 I_{m_1}\cdots\oplus b_k
I_{m_k}$, where $a_1,b_1,\dots,b_k$ are distinct. By Theorem \ref{k} $C_{(X-a_1)^2}\oplus a_1 I_{m}\sim g(C_f)$,
where $g=X^{m+1}+a_1$, $f=X^{m+2}$. Let $0,c_1,\dots,c_k\in F$ be distinct. Then
$b_i I_{m_i}=g_i(C_{f_i})$ for $f_i=(X-c_i)^{m_i}$, $g_i=b_i$, $1\leq i\leq k$, so
Lemma \ref{chin} applies.  $\quad\blacksquare$

\begin{thm}
\label{rojo}
Let $F$ be a field that has at least $n-2$ elements and let $A\in M_n(F)$.
Then $A\sim g(C_f)$ for suitable polynomials $f,g\in F[X]$,
with $f$ monic of degree~$n$.
\end{thm}

\noindent{\it Proof.} We may assume that we are not in any of situations described in Lemmas \ref{xr} and \ref{xr2}.
We may also assume that $A$ itself is not similar to a companion matrix.

Let $P=\{p_1,\dots,p_k\}$ be the set of monic irreducible factors of the minimal polynomial of~$A$.
In view of our assumptions we have
$\ell(p)\leq n-2$ for every $p\in  P$. If
$\ell(p)< n-2$ for all $p\in P$, or all $p\in P$ are linear (in which case their $F^+$-stabilizer is trivial), then
Theorem \ref{main} applies. Suppose then that $\ell(p)=n-2$ for some $p\in P$, and some $q\in P$ is not linear.
Then either $n=4$, $p=q$ has degree 2, and $A\sim C_p\oplus C_p$, in which case Proposition \ref{z2} applies, or else $p$ is linear, $q$
has degree 2 and the elementary divisors of $A$ are all linear except for $q$, in which case Theorem \ref{main} also applies.
$\quad\blacksquare$

\section{The case when all $p$-elementary divisors of $A$ are equal}

Let $A\in M_n(F)$ and let $P=\{p_1,\cdots,p_m\}$ be
the set of monic irreducible polynomials of the minimal polynomial of $A$.
Recall that $A$ is semisimple if its minimal
polynomial is square free. This means that for $p\in P$
the $p$-elementary divisors of $A$ are $p,\dots,p$. We will say that $A$ is {\em homogeneous}
if, given any $p\in P$, there exist positive integers $i=i(p)$ and $k=k(p)$ such that the $p$-elementary divisors of $A$ are $p^i,\dots,p^i$, with $k$ repetitions.

R. Guralnick (private communication) proved that every semisimple matrix is of polynomial type.
Using his idea, we have extended this result to all homogeneous matrices.

\begin{lem}
\label{potq} Let $m\geq 1$ and let $f=a_0+a_1X+\cdots+a_{m-1}X^{m-1}+X^m\in
F[X]$, where $a_0\neq 0$. Given $k\geq 1$, let $B\in M_{km}(F)$ be a matrix made up
of $k^2$ blocks of size $m\times m$, where the diagonal and first superdiagonal blocks are equal to $C_f$,
and all other blocks are 0. Let $e_1,\dots,e_{km}$ be the canonical basis of the column space $V=F^{km}$.
Then the minimal polynomial of $e_{(k-1)m+1}$ relative to $B$ is $f^k$. In particular, $B$ is cyclic
with minimal polynomial $f^k$.
\end{lem}

\noindent{\it Proof.} By induction on $k$. The case $k=1$ is clear. Suppose $k>1$ and the result is true
for $k-1$. We have $V=U\oplus W$, where $U$ (resp. $W$) is the span of the first $(k-1)\times m$ (resp. last~$m$)
of the canonical vectors $e_1,\dots,e_{km}$. Note that $U$ is $B$-invariant. Let $u_1,\dots,u_m$ be the last $m$ of the stated spanning vectors of $U$ and let $v_1,\dots,v_m$ be the stated spanning vectors of~$W$. In particular,
$v_1=e_{(k-1)m+1}$.

Let $g$ be the minimal polynomial of $v_1$ relative to $B$. Clearly the conductor of $v_1$ into $U$ relative to $B$ is $f$, so $g=hf$ for some $h\in F[X]$. On the other hand, we have
$$
Bv_1=u_2+v_2,\quad B^2v_1=u_3+v_3,\dots, B^{m-1}v_1=u_m+v_m,
$$
$$
B^m v_1=-(a_0 u_1+a_1 u_2+\cdots+a_{m-1}u_m)-(a_0 v_1+a_1 v_2+\cdots+a_{m-1}v_m),
$$
so
$$
f(B)v_1=(B^m+a_{m-1}B^{m-1}+\cdots+a_1 B+a_0 I)v_1=-a_0u_1.
$$
It follows that
$$
0=g(B)v_1=h(B)f(B)v_1=-a_0h(B)u_1.
$$
Since $a_0\neq 0$ we infer $h(B)u_1=0$. It follows by inductive hypothesis that $f^{k-1}|h$.
But these are monic polynomials of the same degree, so $f^{k-1}=h$, whence $g=f^k$. $\quad\blacksquare$

\begin{thm}
\label{homog}
Every homogeneous -and so every semisimple- matrix is of polynomial type.
\end{thm}

\noindent{\it Proof.} Let $A\in M_n(F)$ and let $P=\{p_1,\cdots,p_m\}$ be
the set of monic irreducible factors of the minimal polynomial of $A$, so that
$A\sim A_{p_1}\oplus\cdots\oplus A_{p_m}$, where for each $p\in P$ there exists $i(p),k(p)\geq 1$ such that
$A_p=C_{p^{i(p)}}\oplus\cdots\oplus C_{p^{i(p)}}$, with $k(p)$ summands. Let $p\in P$. Suppose first $p(0)\neq 0$. Let $B_p$
be the matrix $B$ of Lemma \ref{potq} constructed upon $f=p^{i(p)}$ and $k=k(p)$. By construction $A_p$ commutes with $B_p$,
which is cyclic with minimal polynomial $p^{i(p)k(p)}$ by Lemma \ref{potq}. Suppose next $p(0)=0$, that is, $p(X)=X$. Now $A_X\sim C_{X^i}\oplus\cdots\oplus C_{X^i}$,
with $k$ repetitions, so by Proposition \ref{z2} we have $A_X\sim g(C_f)$, for $g=X^k$ and $f=X^{ki}$.
In particular, $A_X$ commutes with a cyclic matrix $B_X$ with minimal polynomial $X^{ki}$.
Since $B_{p_1},\dots,B_{p_m}$ are cyclic with pairwise relatively prime minimal polynomials, it follows
that $B_{p_1}\oplus\cdots\oplus B_{p_m}$ is also cyclic.  As this matrix commutes with $A_{p_1}\oplus\cdots\oplus A_{p_m}$, the result follows
from Lemma \ref{comiu}.$\quad\blacksquare$

\section{Examples}
\label{cou}

Let $F$ be an arbitrary finite field, say $F=\{a_1,\dots,a_m\}$,
where $a_1=0$. Let $A\in M_n(F)$ be the matrix of size $n=m+3$
equal to the direct sum of the following $t=m+1$ Jordan blocks:
$$
 A=J_3(0)\oplus J_1(0)\oplus J_1(a_2)\oplus\cdots\oplus
J_1(a_m).
$$
We claim that $A$ is not of polynomial type. Suppose, by way of contradiction, that
$A\sim g(C_f)$ for some $f,g\in F[X]$, with $f$ monic of degree
$m+3$. For $1\leq i\leq m$ let $h_i=\gcd(g-a_i,f)$. It follows from Theorem \ref{k} that $h_i$ is linear for $i>1$ and of
degree 2 if $i=1$. Note that $h_1,\dots,h_m$ are pairwise
relatively prime, since so are $g-a_1,\dots,g-a_m$. Let $h=h_1$
and set $u_i=\mathrm{deg}\; gcd (g^i,f)$, $i\geq 1$. By
(\ref{pat}) we have $u_i=\mathrm{deg}\; gcd (h^i,f)$ for all
$i\geq 1$. By Theorem \ref{k}
$$
u_1=2,\; u_2=3,\; u_3=4.
$$
Since $h$ has degree 2 we see from $u_2=3$ that $h$ cannot be irreducible. Then
either $h=r^2$ or $h=rs$, where $r,s\in F[X]$ are distinct, linear and monic. If $h=r^2$ then $u_2=3$
forces $r$ to have multiplicity 3 in $f$, against $u_3=4$. If $h=rs$ then the $m+1$ distinct
monic linear polynomials $r,s,h_2,\dots,h_m$ of $F[X]$ are factors of $f$, which is absurd as $|F|=m$.
This proves our claim.

This shows that over any given finite field there is a matrix that is not of polynomial type.
Since, $|F|<n-2$ and $|F|<t$, our example shows that the
bounds given in Theorems \ref{thor} and~\ref{rojo} are best
possible. Moreover, our example also proves that Lemmas \ref{xr} and \ref{xr2} cannot be pushed any further
in the same direction. Furthermore, the presence of the $X$-elementary divisors $X^3,X$ in $A$ shows
that our conditions on Theorem \ref{homog} are well placed. In this regard, we believe that even the following
might be true: if $p\in F[X]$ is any monic irreducible polynomial then
there exists $A\in M_n(F)$, all whose elementary divisors are powers of $p$,
that is not of polynomial type (see Section \ref{nilo} for the case $p=X$).

R. Guralnick asks whether a semisimple matrix commutes with a cyclic {\em semisimple} matrix.
In the infinite case, or under the weaker hypotheses of Theorem \ref{main}, we see that
any $A\in M_n(F)$ will commute with a cyclic matrix $B$ whose elementary divisors have the same
exponents as those of $A$. In particular, $B$ is semisimple if so is $A$. We do not know the answer
in general. In this regard, consider the following statement.

(S) Let $F$ be a field and let
$P=\{p_1,\dots,p_m\}$ be a set of monic irreducible polynomials in $F[X]$. Let $d_1,\dots,d_m$
be a list of positive integers (with possible repetitions). Then there exist polynomials
$g_1,\dots,g_m\in F[X]$ such that each $g_i$ has degree $d_i$, each $p_i(g_i(X))$ is multiplicity free, and the polynomials $p_1(g_1(X)),\dots,p_m(g_m(X))$ are pairwise relatively prime.

Note that {if} (S) is true when $F$ is a finite field then every semisimple matrix commutes with a cyclic semisimple matrix by Proposition \ref{z2} and Lemma \ref{chin}. For instance,

\noindent $\bullet$ $F=F_2$, $p=X^2+X+1$ and $A=C_p\oplus C_p$. Take $g=X^2+X$. Then $p(g(X))=X^4+X+1$
is irreducible and $A\sim g(C_{p(g(X))})$ by Proposition \ref{z2}.

\noindent $\bullet$ $F=F_3$, $p=X^2+1$ and $A=C_p\oplus C_p\oplus C_p\oplus C_p$. Take $g=X^4$.
Then $p(g(X))=X^8+1$ is multiplicity free (look at its derivative)
and $A\sim g(C_{p(g(X))})$ by Proposition \ref{z2}.

\noindent $\bullet$ $F=F_4$, $p=X^2+X+\a$, where $\a^2+\a+1=0$, $q=X^2+X+(\a+1)$, $A=C_p\oplus C_p\oplus C_q$
and $B=C_p\oplus C_p\oplus C_q\oplus C_q$.  Take $g=X^2+X$. Then $p(g(X))=X^4+X+\a$, which is multiplicity free
and not a multiple of $q$, and $q(g(X))=X^4+X+(\a+1)$ is multiplicity free and relatively
prime to $p(g(X))$. Hence $A$ (resp. $B$) commutes with a cyclic and semisimple matrix with
minimal polynomial $p(g(X))q(X)$ (resp. $p(g(X))q(g(X))$).

\begin{exa}\label{gol}{\rm Let $F=F_p$, where $p$ is a prime, and consider an $F^+$-invariant irreducible polynomial of degree $p$, which must necessarily be of the form
$q(X)=X^p-X-a$, where $0\neq a\in F$.
We claim that $A=C_{q^2}\oplus C_q$ is of polynomial type.

We consider first the case $p$ odd. Make the initial assumption that $a=1$.
Let $g(X)=X^2$ and consider the polynomial $h(X)=q(g(X))=X^{2p}-X^2-1$. We claim that $h$ is the
product of 2 irreducible relatively prime polynomials $r$ and $s$ of degree $p$.
Since $h'(X)=-2\neq 0$, it is clear
that $h$ has no repeated irreducible factors. Now $q(X)=(X-\a)(X-\a^p)\cdots (X-\a^{p^{p-1}})$ for a root $\a$
in some extension of $F_p$. Thus
$$
h(X)=q(X^2)=(X^2-\a)(X^2-\a^p)\cdots (X^2-\a^{p^{p-1}}).
$$
Clearly
$$
1=\a\a^p\cdots\a^{p^{p-1}}.
$$
Since $1+p+p^2+\cdots+p^{p-1}$ is odd, $\a$ has a square root, say $\b$, in the multiplicative group of the field $F_p[\a]$. Thus
$$
h(X)=(X-\b)(X-\b^p)\cdots (X-\b^{p^{p-1}})\times (X+\b)(X+\b^p)\cdots (X+\b^{p^{p-1}})=r(X)\times s(X)
$$
The polynomials $r,s$ are clearly invariant under the Frobenius automorphism $x\mapsto x^p$, so
they have coefficients in $F_p$. They are irreducible, since the degree of $\b$ and $-\b$ is $p$
(the degree cannot be 1 since $\beta^2=\a$, and there is simply no other option for this degree).
This establishes the claim.

We now let $f(X)=r^2(X)s(X)$. Then $f(X)$ divides $q^2(g(X))=h(X)^2=r(X)^2 s(X)^2$ but not $q(g(X))=r(X)s(X)$.
Since $q(X)$ is irreducible, the minimal polynomial of $g(C_f)$ is $q^2$ and the remaining invariant factor has no
other option than being $q$. Thus $A\sim g(C_f)$. Now in the general case $q=X^p-X-a$, we claim that the same
choice of $f$ works with $g(X)=aX^2$. Indeed, $q(g(X))=aX^{2p}-aX^2-a=a(X^{2p}-X-1)$. Thus, by the previous case,
$f$ divides $q^2(g(X))$ but not $q(g(X))$, as required.

As an illustration, when $p=3$ and $q(X)=X^3-X-1$ we have $h(X)=X^6-X^2-1=(X^3-X^2-X-1)(X^3+X^2-X+1)=r(X)s(X)$.
For $q(X)=X^3-X+1$ we take $g(X)=-X^2$.

In the case $p=2$ we have $q(X)=X^2+X+1$ and we just take $f(X)=q(X)^3$ and $g(X)=X^2$. Then $f$ divides $q^2(g(X))=q^4(X)$
but not $q(g(X))=q^2(X)$, as required.}
\end{exa}

\section{The case of a nilpotent matrix}
\label{nilo}

Here we focus on a matrix $A\in M_n(F)$ that is equal to the direct
sum of Jordan blocks of various sizes corresponding to the same
eigenvalue. We wish to see if $A$ is of polynomial type. In view of our next result
we may assume that this eigenvalue is 0, i.e. $A$ is
nilpotent.

\begin{lem}
\label{y3} Suppose that $A$ is similar to $g(C_f)$. Then $A+aI$ is
similar to $(g+a)(C_f)$.
\end{lem}

\noindent{\it Proof.} We have $YAY^{-1}=g(C_f)$, for
$Y\in\GL_n(F)$, so $Y(A+aI)Y^{-1}=(g+a)(C_f)$.$\quad\blacksquare$

\medskip

Assume henceforth that $A$ is nilpotent and write $A=1^{c(1)}|\cdots |n^{c(n)}$, $c(i)\geq 0$, to mean that
$A$ is the direct sum of $c(i)$ 0-Jordan blocks $J_i$ of size $i$, where $1\cdot
c(1)+2\cdot c(2)+\cdots+n\cdot c(n)=n$.

There may be repetitions amongst the non-zero $c(1),\dots,c(n)$ and this will play an important role. To account
for this, we set $r(i)=|\{j\,|\, c(j)=c(i)\}|$ for all $1\leq i\leq n$ such that $c(i)\neq 0$.

Suppose that $F[X]$ has at least $r(i)$ monic irreducible polynomials of degree $c(i)$, for all $1\leq i\leq n$
such that $c(i)\neq 0$. We may then select monic polynomials $q_1,\dots,q_n\in F[X]$ such that $\mathrm{deg}(q_i)=c(i)$, with $q_i$
irreducible if $c(i)\neq 0$, and all polynomials $q_1,\dots,q_n$ different from 1 are distinct from each other. Using
Theorem \ref{k} it is easy to see that
\begin{equation}
\label{irsim}
A\sim g(C_f),\;\text{for }g=q_1q_2\cdots q_n\text{ and }f=q_1^1 q_2^2\cdots q_n^n.
\end{equation}
Observe how (\ref{irsim}) agrees with Proposition \ref{z2} when all but one of the indices $j$ satisfy $c(j)=0$,
i.e. when $A=i^c$ (from now on we agree to omit $j$ if and only if $c(j)=0$, as well as $c(j)$ if and only if $c(j)=1$).

Let us concentrate on the case when $F$ is a finite field (otherwise Theorem \ref{rojo} settles the matter).
For each $d\geq 1$, let $N(d)$ be the number of monic irreducible polynomials of degree $d$ in $F[X]$.
There is a formula (see [J2]) due to Gauss for $N(d)$. Using Gauss' formula we may verify whether or not
\begin{equation}
\label{gauss}
N(c(i))\geq r(i),\text{ when } c(i)\neq 0.
\end{equation}
However, as seen below, condition (\ref{gauss}) is not necessary for $A$ to be of polynomial type, as one may establish this fact
using other polynomials $f,g$. It is however necessary in some cases, and this allows us to provide examples of nilpotent matrices
that are not of polynomial type.

Suppose first $A=1|2|i$, i.e. $A=J_1\oplus J_2\oplus
J_i$. By Theorem \ref{k} $A\sim g(C_f)$ with $g=r^2s$,
$f=r^3s^i$ and $r\neq s\in F[X]$ linear. Here $r(1)=3$ and $F_2[X]$ has only 2 linear polynomials, but
we can still show that $A$ is of polynomial type using a different method than (\ref{irsim}).

Suppose next $A=1|3|4$, i.e. $A=J_1\oplus J_3\oplus
J_4\in M_8(F)$. By Theorem \ref{k} $A\sim g(C_f)$ with $g=r^2s$,
$f=r^7s$ and $r\neq s\in F[X]$ linear. The same comments made above apply.

This cannot be continued. Indeed, suppose next $F=F_2$ and $A=1|3|5\in M_9(F)$. In this case $N(1)=2<3=r(1)$. We claim that $A$ is not of polynomial type. To
see the claim suppose, by way of contradiction, that $A\sim
g(C_f)$ for some $f,g\in F[X]$, with $f$ monic of degree~$9$. Let
$h=\gcd(g,f)$ and set $u_i=\mathrm{deg}\; \gcd (g^i,f)$, $i\geq 1$.
By (\ref{pat}) we have $u_i=\mathrm{deg}\; \gcd (h^i,f)$ for all
$i\geq 1$. By Theorem \ref{k}
$$
u_1=3,\; u_2=5,\; u_3=7,\; u_4=8,\; u_5=9
$$
Thus $h$ has degree 3. From $u_2=5$ we see that $h$ cannot be
irreducible. Suppose $h=rs$, with $r,s\in F[X]$ irreducible of
degrees 2 and 1, respectively. Then $u_2=5$ forces $s$ to have
multiplicity 1 in $f$, but then $u_4=8$ becomes impossible.
Suppose next $h=r^2 s$, with $r,s\in F[X]$ distinct and linear.
Then $u_2=5$ forces $s$ to have multiplicity 1 or $r$ multiplicity
3 in $f$. In the first case $u_4=8$ forces $r$ to have
multiplicity 7 in $f$, which contradicts $u_5=9$. In the second case $u_3=7$ becomes
impossible. This proves the claim.

Using the preceding ideas we easily see that every nilpotent matrix of size
$n\leq 8$ over any field is of polynomial type. The example of $A=1|3|5$ over $F_2$ shows that this
cannot be extended to $n=9$.

Using exactly the same method we also verified that if $F=F_3$ and $A=1|4|7|10$ then $A$ is
not of polynomial type. This seems to indicate that if $F$ has size $m$ and
$A=1|m+1|\cdots|m^2+1$ then $A$ is not of polynomial type.

As seen above, while (\ref{irsim}) gives a pleasant formula, it is not sufficiently general
to deal with all cases in which $A$ is of polynomial type. This can be achieved, at least
in principle, using our work from section \ref{jam}, as hinted in our treatment of the case $A=1|3|5$.
Indeed, Claim \ref{calmp} shows that the similarity classes of
nilpotent matrices $A$ in $M_n(F)$ of fixed nilpotency class $t$,
with $1\leq t\leq n$, are in bijective correspondence with the
collection of all strictly increasing sequences of natural numbers
$s_1<s_2\dots< s_t=n$. If $N(A)$ is the nullity of $A$, the
correspondence is
$$
[A]\mapsto s_1=N(A)<s_2=N(A^2)<\cdots<s_t=N(A^t).
$$
By Theorem \ref{k} such $A$ will be of polynomial type if and only if there exist
$f,g\in F[X]$, with $f$ monic of degree~$n$, such that
$$
\mathrm{deg}\;\gcd(g^i,f)=s_i,\quad 1\leq i\leq t.
$$
By (\ref{pat}) we may assume that $g|f$. Let $f=r_1^{a_1}\cdots
r_k^{a_k}$ and $g=r_1^{b_1}\cdots r_k^{b_k}$ be the prime
decompositions of $f$ and $g$, where $1\leq a_1$, $b_i\leq a_i$
and $d_i=\mathrm{deg}(r_i)$ for all $1\leq i\leq k$. It is
necessary to require that $1\leq b_i$ for all $1\leq i\leq k$
since $g(C_f)^t=0$ implies $f|g^t$.

Thus, $A$ will be of polynomial type if and only if there exist
$k\geq 1$ distinct monic irreducible polynomials $r_1,\dots,r_k$
in $F[X]$ of degrees $d_1,\dots,d_k$ and integers $1\leq b_i\leq
a_i$, for all $1\leq i\leq k$, such that
$$
\mathrm{min}\{i
b_1,a_1\}d_1+\cdots+\mathrm{min}\{ib_k,a_k\}d_k=s_i,\quad 1\leq
i\leq t.
$$
The first of these equations says that $\mathrm{deg}(g)=s_1$ and
and the last that $f|g^t$.

\section{Further consequences}

We present applications of our work in section \ref{jam} which were postponed to deal
with  matrices of polynomial type.

Our first application, to linear algebra itself, follows easily from our work in section \ref{jam}.

\begin{thm}
\label{csd}
 Let $g\in F[X]$ and let $A\in M_n(F)$ have minimal polynomial $f$.
 Let $P$ be the set of minimal polynomials of $g$ modulo the monic
irreducible factors of $f$. Then $g(A)$ is semisimple if and only if $\gcd(p(g(X)),f)=\gcd(p^2(g(X)),f)$ for all $p\in P$. Moreover,
$g(A)$ is cyclic if and only if $A$ is cyclic and
$\mathrm{deg}\;\gcd(p(g(X)),f)=\mathrm{deg}\; p$ for all $p\in P$. Furthermore, all eigenvalues of $g(A)$ are in $F$ if and only if $P=\{X-a_1,\dots, X-a_m\}$
for some $a_1,\dots,a_m\in F$, in which case these are the eigenvalues of $g(A)$, with $g(A)$ diagonalizable if and only if  $\gcd(g-a_i,f)=\gcd((g-a_i)^2,f)$  for all $1\leq i\leq m$.
\end{thm}

We next give an application to simple algebraic field extensions. Let $F$ be a field, let $f\in F[X]$
be monic of degree $n$ and let $\a$ be algebraic over $F$ with minimal polynomial $f$. Then $K=F[\a]$
is a field with $F$-basis $B=\{1,\a,\dots,\a^{n-1}\}$. Given $\b\in K$ one often needs to find the minimal
polynomial, trace, norm, or inverse of $\beta$ (in this last case if $\beta\neq 0$). There are several
ways of doing this. One way is to look at the regular representation $\ell:K\to \mathrm{End}_F(K)$,
where $\ell_\beta$ is the map ``multiplication by $\beta$" and find the previous information about $\ell_\b$,
using the matrix $M_B(\ell_\b)$ of $\ell_\b$ relative to $B$. We wish to find a closed formula
for $M_B(\ell_\b)$ that can be of use for all these purposes. We have $\beta=g(\a)$ for a unique $g\in F_n[X]$,
and it is clear that $M_B(\ell_\a)=C_f$. Therefore, invoking Lemma 3.1 of \cite{GS2}, we obtain
\begin{equation}
\label{uqa}
M_B(\ell_\beta)=M_B(\ell_{g(\a)})=M_B(g(\ell_\a))=g(M_B(\ell_\a))=g(C_f)=([g]\,C_f[g]\dots \,C_f^{n-1}[g]).
\end{equation}
Thus $M_B(\ell_\beta)$ can be computed by applying the successive powers of $C_f$ to $[g]$. It is also
clear that (\ref{uqa}) remains valid if we replace $F$ by an arbitrary commutative ring with identity $R$,
as long as $f\in R[X]$ satisfies $f(\a)=0$ and $g(\a)\neq 0$
for all $g\in R_n[X]$.

The matrix $g(C_f)=([g]\,C_f[g]\dots \,C_f^{n-1}[g])$ can be used to answer
any questions about $\beta=g(\a)$. In the field case, where the minimal polynomial $p$ of $\beta$ is irreducible,
our work from Section \ref{jam} then says that $p$ is the minimal polynomial of $g$ modulo $f$
(this is how one finds $p$ in practice, perhaps without noticing). The trace and norm of $\beta$ are those of $g(C_f)$. If $\beta\neq 0$ then $\beta^{-1}=[h]$, where $h\in F_n[X]$ and $[h]$ is the first column of $g(C_f)^{-1}$.
Of course, $g(C_f)^{-1}=h(C_f)$, and $h$ can also be computed using Euclid's algorithm. If we actually
compute $p$ first, as the minimal polynomial of $g$ modulo $f$, we can find $\beta^{-1}$ in a third way,
namely by using the relation $p(\beta)=0$ and isolating the independent term.

\medskip

It turns out that Lemma \ref{uno} and some of its consequences hold in the context of rings.
For the remainder of this section we assume that $R$ is an integral domain with field of fractions $F$.
Recall that $R$ is said to be integrally closed if any root in $F$ of a monic polynomial in $R[X]$ actually lies in $R$.

\begin{lem}
\label{ic}
 Suppose that $R$ is integrally closed. Let $f,g\in F[X]$, with $f$ monic, and
let $d=\gcd(f,g)\in F[X]$. Then $d\in R[X]$.
\end{lem}

\noindent{\it Proof.} Let $K$ be an algebraic closure of $F$. Then
$d=(X-a_1)\cdots (X-a_m)$ for some $a_i\in K$. The $a_i$ are roots
of $f$ and hence they are integral over $R$. It follows that the
coefficients of $d$ are integral over $R$. These coefficients are
in $F$ and $R$ is integrally closed, so $d\in R[X]$.$\quad\blacksquare$

\medskip

We now extend all of the notation used in Lemma \ref{uno} to $R$. Assume
that $z=\gcd(y(g(X),f)$, as calculated in $F[X]$, actually belongs to $R[X]$
(this is guaranteed if $R$ is integrally closed by Lemma \ref{ic}).
Then Lemma \ref{uno} is true in this generality. Thus, $[p]\in V(y)$ if and
only if $h|p$, and $V(y)$ is a free $R$-module of rank $d=\mathrm{deg}\; z$ with $R$-basis $[h],[X h],\dots, [X^{d-1} h]$. The same proof works.

The following is a consequence of the preceding remarks and Lemma \ref{uno} applied to the case $g=X$ and $y(X)=(X-\a)^i$.

\begin{cor}
\label{equis} Let $\a$ be a root of
$f$ in $R$ of multiplicity at least $i$, denote the corresponding generalized eigenspace of $C_f$ by  $E=\{v\in V\,|\,
(C_f-\a I)^i v=0\}$,
and set $h=f/(X-\a)^i$. Then $E$ is a free $R$-module of rank $i$
spanned by $[h], [h X],\dots, [X^{i-1} h]$.
\end{cor}

The following is a consequence of Corollary \ref{equis} applied to the case $i=1$.

\begin{cor}
\label{equw}
 Let $\a$ be a root of $f$ in $R$, let $E$ be corresponding eigenspace of $C_f$
and set $h=f/(X-\a)$. Then $E$ is a free $R$-module of rank 1 spanned by $[h]$.
\end{cor}

Note that the field cases of Corollaries \ref{equis} and
\ref{equw} appear in different but equivalent forms in \cite{B}, \cite{Gi}
and \cite{R}.

The following is a consequence of Lemma \ref{uno} applied to the
case $y(X)=X$ and arbitrary $g\in R[X]$.

\begin{cor}\label{upa}
Let $N$ be the nullspace of $g(C_f)$.
Suppose that $z=\gcd(g,f)\in F[X]$ actually belongs to $R[X]$ (this is guaranteed if $R$ is integrally closed) and let $d$ be its degree. Set $h=f/z$. Let $p\in R_n[X]$. Then $[p]\in N$ if and only if $h|p$.
Moreover, $N$ is a free $R$-module of rank $d=\mathrm{deg}\; z$ and, if $d>0$, then $[h],[X h],\dots, [X^{d-1} h]$
is an $R$-basis of $N$.
\end{cor}


\section{Companion matrices and subalgebras of $M_n(R)$}
\label{ante}

In this section $R$ is an integral domain with field of fractions $F$. Let $f,g\in R[X]$
be monic polynomials of degree $n$, with respective companion matrices $C,D\in M_n(R)$. We further let $d=\gcd(f,g)\in F[X]$,
whose degree will be denoted by $m$. Note that $d\in R[X]$ provided $R$ is integrally closed, by Lemma \ref{ic}.
We also let $h=f/d\in F[X]$.

We let $R\langle C,D\rangle$ stand for the subalgebra of $M_n(R)$ generated by $C$ and
$D$, and $R[C,D]$ for the $R$-span of all matrices $C^iD^j$,
where $0\leq i,j<n$. It is shown in Corollary 6.3 of \cite{GS2} that $R\langle C,D\rangle=R[C,D]$. Moreover, Theorem 9.2 of
\cite{GS2} shows that if $R$ is a unique factorization domain then $R\langle C,D\rangle$ is necessarily  a free
$R$-module of rank $n+(n-m)(n-1)$, and Theorem 9.3 of
\cite{GS2} gives a presentation for
$R\langle C,D\rangle$ in terms of $C,D$. As indicated in the Introduction, we wish to extend these results
to the more general setting of
integrally closed domains as well as furnishing a counterexample to the freeness of $R\langle C,D\rangle$ when $R=\Z[\sqrt{5}]$.

\begin{thm}
\label{piu}
 The $n+(n-m)(n-1)$ elements
$$
I,C,...,C^{n-1},D,CD,...,C^{n-m-1}D,...,D^{n-1},CD^{n-1},...,C^{n-m-1}D^{n-1}
$$
of $R\langle C,D\rangle$ are linearly independent over $F$, and hence over $R$.
\end{thm}

\noindent{\it Proof.} Suppose that $p_0\in F_n[X]$ and $p_1,...,p_{n-1}\in F_{n-m}[X]$ satisfy
\begin{equation}
\label{x}
p_0(C)I+p_1(C)D+\cdots+p_{n-1}(C)D^{n-1}=0.
\end{equation}
We wish to show that $p_0,p_1,...,p_{n-1}=0$. Since $\mathrm{deg}\; p_0<\mathrm{deg}\; f$ and
$\mathrm{deg}\; p_i<\mathrm{deg}\; h$ for $i\geq 1$, it suffices to show that $f|p_0$ and $h|p_1,...,h|p_{n-1}$. Since
$$
p_1(C)D+\cdots+p_{n-1}(C)D^{n-1}=-p_0(C),
$$
it follows that $p_1(C)D+\cdots+p_{n-1}(C)D^{n-1}$ commutes with $C$. Thus
$$
p_1(C)(CD-DC)+p_1(C)(CD^2-DC^2)+\cdots+p_{n-1}(C)(CD^{n-1}-D^{n-1}C)=0.
$$
Now
$$
(CD-DC)e_1=\cdots=(CD^{n-2}-D^{n-2}C)e_1=0,
$$
so for $s=f-g\in R_n[X]$ we have
$$
0=p_{n-1}(C)(CD^{n-1}-D^{n-1}C)e_1=p_{n-1}(C)[s].
$$

By Lemma 3.2 of \cite{GS2}, $f$ divides $p_{n-1}s$ and hence
$p_{n-1}g$. It follows that $h$ divides $p_{n-1}$, so
$p_{n-1}=0$. Proceeding like this with $e_2,...,e_{n-1}$ we see
that $p_{n-2}=p_{n-3}=...=p_1=0$. Going back to (\ref{x}) shows
that $f|p_0$, so $p_0=0$.$\quad\blacksquare$

\begin{cor}
\label{pqw}
 The $n+(n-m)(n-1)$ matrices
$$
I,C,...,C^{n-1},D,CD,...,C^{n-m-1}D,...,D^{n-1},CD^{n-1},...,C^{n-m-1}D^{n-1}
$$
span $R\langle C,D\rangle$ if and only if $d\in R[X]$. In particular, $R\langle C,D\rangle$ is a free $R$-module of rank $n+(n-m)(n-1)$ provided $d\in R[X]$, and in particular if $R$ is integrally closed.
\end{cor}

\noindent{\it Proof.} Recall first of all that $R\langle C,D\rangle=R[C,D]$.

Suppose first $d\in R[X]$. It follows that $h$ also belongs to
$R[X]$. By Lemma 9.1 of \cite{GS2}, if $p\in R[X]$ then dividing
$p$ by $h$ we may write $p(C)D$ in the form $q(C)+z(C)D$, where
$q\in R[X]$ and $z\in R_{n-m}[X]$, so the listed matrices span
$R\langle C,D\rangle$.

If $d\notin R[X]$ then $h\notin R[X]$. The linear independence of
the listed matrices and Lemma~9.1 of \cite{GS2} show that
$C^{n-m}D$ is an $F$-linear combination of the listed matrices but
not an $R$-linear combination of them.$\quad\blacksquare$

\begin{cor}
\label{hg} Let $R$ be an integral domain. If $R\langle C,D\rangle$  is an $R$-free module then its rank must be $n+(n-m)(n-1)$.
\end{cor}

\noindent{\it Proof.} An $R$-basis of $R\langle C,D\rangle$ is an $F$-basis of $F\langle C,D\rangle$, which has dimension $n+(n-m)(n-1)$ by Corollary \ref{pqw}. $\quad\blacksquare$

\begin{thm} Suppose $d\in R[X]$. Let the polynomials $P_1,\dots,P_{n-1}\in R[X,Y]$
be defined as in Theorem 7.1 of \cite{GS2} or, more generally, be
arbitrary while satisfying
$$D^jC = P_j(C,D),\quad j = 1, . . . , n-1.$$
 Then the algebra $R\langle
C,D\rangle$ has presentation
$$
\langle X,Y\,|\, f(X)=0,\, g(Y)=0,\, h(X)(X-Y)=0,\,
Y^jX=P_j(X,Y),\quad j=1,\dots,n\rangle.
$$
\end{thm}

\noindent{\it Proof.} Write $\Omega:R\langle X,Y\rangle\to
R\langle C,D\rangle$ for the natural $R$-algebra epimorphism that
sends $X$ to $C$ and $Y$ to $D$. Let $K$ be the kernel of
$\Omega$. Set $S=R\langle X,Y\rangle/K$ and let $A$ and $B$ be the
images of $X$ and $Y$ in $S$. Then $g(B)=0$ and $B^jA=P_j(A,B)$,
$j=1,\dots,n-1$. It follows from  Lemma~6.1 of \cite{GS2} that
$S=R[A,B]$. From $f(A)=0$ and $g(B)=0$ it follows that $S$ is
$R$-spanned by $A^iB^j$, $0\leq i,j\leq n-1$. Thus the relation
$h(A)(A-B)=0$ allows $R[A,B]$ to be spanned by the reduced list of
$n+(n-m)(n-1)$ matrices:
$$
I,A,\dots,A^{n-1},
B,AB,\dots,A^{n-m-1}B,\dots,B^{n-1},AB^{n-1},\dots,A^{n-m-1}B^{n-1}.
$$

Let $\Delta:S\to R\langle C,D\rangle$ be the map induced by $\Omega$.
We wish to show that $\Delta$ is injective. Let $t\in\ker \Delta$. Then
$t$ is a linear combination of the above listed matrices. Their images under
$\Delta$ are linearly independent by Theorem \ref{piu}, so $t=0$. $\quad \blacksquare$

\medskip

Our next example shows that $R\langle C,D\rangle=R[C,D]$ need not be a free $R$-module.

\begin{exa} Let
$R=\Z[\go]$, with field of fractions $F=\Q[\go]$. Let $n=2$ and set
$$
g=(X-(1-\go)/2)(X-(1+\go)/2)=X^2-X-1\in R[X],
$$
$$ f=(X-(1-\go)/2)(X-(5+\go)/2)=X^2-3X-\go\in R[X].
$$
Then the subalgebra $R\langle C,D\rangle$ of $M_2(R)$ is not a free $R$-module.
\end{exa}

\noindent{\it Proof.} This will proceed in a number of steps.

\medskip

\indent (S1) Let $L$ be the ideal of $R$ generated by $2$ and $1-\go$ and let $J$ be the subset of $R$ formed
by all $a+b\go$ such that $a,b\in\Z$ have the same parity, i.e. $a\equiv b\mod 2$. We claim
that $L=J$ and that $R/L\cong \Z_2$ is a field with 2 elements.

First of all $J$ is clearly a subgroup of $R$ that is stable under multiplication by $\go$.
It follows that $J$ is an ideal of $R$. Since $2,1-\go\in J$, we conclude that $L\subseteq J$.
On the other hand, all $2a+2b\go$ are clearly in $L$ as well as all $2a+1+(2b-1)\go=2(a+b\go)+1-\go$.
Thus $J\subseteq L$, so $L=J$. Finally, $\go\equiv 1\mod L$, so $a+b\go\equiv a+b\equiv 0,1\mod L$, as required.

\medskip

\indent (S2) $R$ is not integrally closed, as $(1-\go)/2\in F$ is a root of the monic polynomial $g\in R[X]$.

\medskip

\indent (S3) The companion matrices of $f$ and $g$ are
$$
C=\left(\begin{matrix} 0 & \go\\ 1 & 3
    \end{matrix}\right),\quad D=\left(\begin{matrix} 0 & 1\\ 1 & 1
    \end{matrix}\right).
$$
Observe that
$$
CD=\left(\begin{matrix} \go & \go\\ 3 & 4
    \end{matrix}\right)
$$
Note that $d=\gcd(f,g)\in F[X]$ is $X-(1-\go)/2$ and $d$ does not belong to $R[X]$. Its degree is $m=1$.
By Theorem \ref{piu}, or by inspection, the matrices $I,C,D$ are linearly independent. By Corollary \ref{hg}, if $R\langle C,D\rangle$ is a free $R$-module its rank must be 3.

\medskip

\indent (S4) Set $M=R[C,D]=R\langle C,D\rangle$ and $K=R/L$, both of which are $R$-modules. Then
$$
H=\mathrm{Hom}_R(M,K)
$$
is an $R$-module via
$$
(rf)(m)=r(f(m)).
$$
Since $L$ annihilates $K$ it also annihilates $H$, which makes $H$ into an $R/L$-modules, i.e. a $K$-vector space.

\medskip

\indent (S5) Here we prove that the dimension of $H$ over $K$ is
4. More explicitly, we show that any $R$-homomorphism $M\to K$ can
be arbitrarily defined on $I,C,D,CD$ (taking a value of 0 or 1 in
$K$) and then extended $R$-linearly to $M$. This is true in spite
of the fact that $CD$ is an $F$-linear combination of $I,C,D$.

Indeed, we know that $M$ is $R$-spanned by $I,C,D,CD$.
Let $x_1,x_2,x_3,x_4\in K$ be any values. Then the necessary and sufficient condition for the existence
of an $R$-homomorphism $f:M\to K$ satisfying
$$f(a_1 I+a_2 C+a_3D+a_4 CD)=a_1x_1+a_2x_2+a_3x_3+a_4x_4,$$
for all $a_i\in R$, is that
$$
a_1 I+a_2 C+a_3D+a_4 CD=0\implies a_1x_1+a_2x_2+a_3x_3+a_4x_4=0.
$$

The curious fact here is that if $a_1 I+a_2 C+a_3D+a_4 CD=0$ then all $a_i\in L$, which clearly
implies $a_1x_1+a_2x_2+a_3x_3+a_4x_4=0$ since $L$ annihilates $K$.

Thus, we are reduced to show that $a_i\in R$ and $a_1 I+a_2 C+a_3D+a_4 CD=0$ implies $a_i\in L$,
i.e. the linear relations between $I,C,D,CD$ all have coefficients in $L$.

To prove this we first work over $F$. We look for all 4-tuples $(a_1,a_2,a_3,a_4)\in F^4$
such that $a_1 I+a_2 C+a_3D+a_4 CD=0$. Since $I,C,D$ are linearly independent and $CD$ is an $F$-linear
combination of $I,C,D$ (by Lemma 9.1 of [GS]) it follows that the space of all these 4-tuples is a subspace
of $F^4$ of dimension 1. Thus all solutions are scalar multiples of a fixed solution. One such solution is
$$
\go I + (1-\go)/2C+(5+\go)/2 D-CD=0.
$$
Let $c\in F$. Scaling the above solution by $c$ yields a solution in $R^4$ if and only if $c=a+b\go\in R$
(look at the coefficient of $CD$) and $c(5+\go)/2\in R$ (look at the coefficients of $C$ and $D$). Now
$$
c(5+\go)/2=(a+b\go)(5+\go)/2=(5(a+b)+(a+5b)\go)/2,
$$
so the last condition can be replaced by $a\equiv b\mod 2$, i.e. $c\in L$. It now only
remains to show that if $c\in L$ then $c\go, c(3-(5+\go)/2),(5+\go)/2,-c\in L$. This reduces
to checking $c(5+\go)/2\in L$, where $c=a+b\go\in L$. If $a=2k,b=2l$ are both even then
$c(5+\go)/2=(k+l\go)(5+\go)\in L$ since $5+\go\in L$. If $a=2k+1,b=2l+1$ are both odd then
$$
c(5+\go)/2=(k+l\go)(5+\go)+(1+\go)(5+\go)/2=(k+l\go)(5+\go)+5+3\go\in L.
$$

\medskip

\indent (S6) $M=R\langle C,D\rangle$ is not a free $R$-module.

Suppose it is free. By (S3) $M$ has rank 3. Then $H$, as defined in (S4), would easily be seen to have dimension 3. But it has dimension 4 by (S5). This completes the proof. $\quad \blacksquare$

\section{Viewing $M_n(R)$ as an $R[X]$-module}

We maintain the notation introduced at the beginning of section
\ref{ante}. Having dealt so extensively with elementary divisors,
it would be unfair not to translate the content of Theorem
\ref{piu} into the language of invariant factors. Note first of
all that the uniqueness part of the invariant factor theorem holds
in much greater generality than its existence counterpart. In
particular it is valid for any commutative ring with identity (as
shown in [K]). We can make $M_n(R)$ into an $R[X]$ module via $C$
as follows:
$$
p(X)\cdot A=p(C)A,\quad p(X)\in R[X], A\in M_n(R).
$$
Clearly $R[C]\subseteq R\langle C,D\rangle=R[C,D]$ are
$R[X]$-submodules of $M_n(R)$, so we can form the quotient module
$M=R\langle C,D\rangle/R[C]$. We further let $D_1=D+R[C]\in
M,\dots, D_{n-1}=D^{n-1}+R[C]\in M$.

\begin{thm} Suppose that $h\in R[X]$ (again, this is automatic if
$R$ is integrally closed). Then the $R[X]$-module $M$ has the
following cyclic decomposition:
$$
M=R[X]\cdot D_1\oplus\cdots R[X]\cdot D_{n-1},
$$
where the annihilating ideals are $R[X]h,\dots,R[X]h$, with $n-1$
repetitions, i.e. the invariant factors of $M$ are $h,\dots,h$,
with $n-1$ repetitions.
\end{thm}

\noindent{\it Proof.} This follows immediately from Theorem
\ref{piu}, its proof, and Lemma 9.1 of
\cite{GS2}.$\quad\blacksquare$

\medskip

Recall from Corollary 5.2 of \cite{GS2} that $M_n(R)=R\langle
C,D\rangle$ if and only if the resultant $R(f,g)$ of $f$ and $g$
is unit in $R$.

\begin{thm}
\label{zv}
 If $R(f,g)$ is a unit in $R$ then the $R[X]$-module $M_n(R)$
has the cyclic decomposition:
$$
M_n(R)=R[X]\cdot C\oplus R[X]\cdot D\cdots\oplus R[X]\cdot
D^{n-1},
$$
where the annihilating ideals are $R[X]f,R[X]f,\dots,R[X]f$, with
$n$ repetitions, i.e. the invariant factors of $M_n(R)$ are
$f,f,\dots,f$, with $n$ repetitions.
\end{thm}

\noindent{\it Proof.} This follows immediately from Theorem
\ref{piu}. $\quad\blacksquare$

Except for the use of $D$, Theorem \ref{zv} is really about the
action of $C$ on $M_n(R)$ and how this yields the invariant
factor $f$ repeated $n$ times. There is no need to use companion
matrices for this: it remains valid in greater generality. For
instance, in the field case, the invariant factors of $M_n(F)$ as
an $F[X]$-module via a given $A\in M_n(F)$ are just the invariant
factors of $A$ repeated $n$ times.

\begin{thm} Let $A\in M_n(R)$ have monic invariant factors $1\neq
f_1|\cdots|f_m$. This means that $A\sim C_{f_1}\oplus\cdots\oplus
C_{f_m}$ or, equivalently, the column space $R^n$, viewed as an
$R[X]$-module via $A$, decomposes as the direct sum of cyclic
submodules with annihilating ideals $R[X]f_1,\dots,R[X]f_m$.

Make $M_n(R)$ into an $F[X]$-module via $A$ by declaring
$p(X)B=p(A)B$, for $p(X)\in R[X]$ and $B\in M_n(R)$. Then $M_n(R)$
has invariant factors
$f_1,\dots,f_1,f_2,\dots,f_2,\dots,f_m,\dots,f_m$, where each
$f_i$ is repeated $n$ times.
\end{thm}

\noindent{\it Proof.} The $F[X]$-module $M_n(R)$ is just the
direct sum of $n$ submodules isomorphic to
$R^n$.$\quad\blacksquare$

\medskip

\noindent{\bf {\large {Acknowledgments}.}} We are indebted to R. Guralnick for his valuable comments,
as well as to the referee who read the manuscript very carefully and made numerous suggestions
that greatly improved the readability of the paper. We also thank D. Farenick, A. Herman and R. Quinlan for their help proofreading the article.
The second author thanks Cecilia, Malena and Federico for their patience during the preparation of the manuscript.


\end{document}